\documentclass[reqno,11pt]{amsart}
\usepackage{amssymb}
\usepackage{amsmath}
\usepackage{mathrsfs}
\usepackage[usenames]{color}
\usepackage{bm}
\usepackage{setspace}
\usepackage{amsaddr}
\numberwithin{equation}{section}
\allowdisplaybreaks

\begin{document}
%\renewcommand{\footnote}{}
%\pagenumbering{roman}
%\thispagestyle{empty}%\quad\newpage
%\thispagestyle{empty}
%\nonstopmode
%**************************************************************************

\newtheorem{theorem}{Theorem}[section] % Nummerierung
\newtheorem{proposition}[theorem]{Proposition}
\newtheorem{corollary}[theorem]{Corollary}
\newtheorem{lemma}[theorem]{Lemma}

\theoremstyle{definition}
\newtheorem{assumption}[theorem]{Assumption}
\newtheorem{definition}[theorem]{Definition}

\theoremstyle{definition} %%{remark}
\newtheorem{remark}[theorem]{Remark}
\newtheorem{remarks}[theorem]{Remarks}
\newtheorem{example}[theorem]{Example}
\newtheorem{examples}[theorem]{Examples}
%**************************************************************************
\newenvironment{pf}%
{\begin{sloppypar}\noindent{\bf Proof.}}%
{\hspace*{\fill}$\square$\vspace{6mm}\end{sloppypar}}
%**************************************************************************
\def\bA{{\bm A}}
\def\bB{{\bm B}}
\def\mE{{\mathbb E}}
\def\mK{{\mathbb K}}
\def\hmE{{\widehat{\mathbb E}}}
\def\mEp{{\mathbb E}_{\phi}}
\def\mFp{{\mathbb F}_{\phi}}
\def\mEpp{{\mathbb E}_{\phi,\mP}}
\def\tmEp{\widetilde{\mathbb E}_{\phi}}
\def\tmEpp{\widetilde{\mathbb E}_{\phi,P}}
\def\mFpp{{\mathbb F}_{\phi,\mP}}
\def\tPhi{\widetilde{\Phi}}
\def\mF{{\mathbb F}}
\def\mG{{\mathbb G}}
\def\mX{{\mathbb X}}
\def\mP{{\mathbb P}}
\def\db{\|}
\def\r{Nr}
\def\R{{\mathbb R}}
\def\N{{\mathbb N}}
\def\C{{\mathbb C}}
\def\Q{{\mathbb Q}}
\def\mP{{\mathbb P}}
\def\Z{{\mathbb Z}}
\def\mH{\mathbb H}
\def\mA{\mathbb A}
\def\mT{\mathbb T}
\def\D{{\mathcal D}}
\def\cB{{\mathcal B}}
\def\E{{\mathcal E}}
\def\cF{{\mathcal F}}
\def\cA{{\mathcal A}}
\def\cH{{\mathcal H}}
\def\G{{\mathcal G}}
\def\B{{\mathcal B}}
\def\I{{\mathcal I}}
\def\M{{\mathcal M}}
\def\O{{\mathcal O}}
\def\S{{\mathcal S}}
\def\cT{{\mathcal T}}
\def\cP{{\mathcal P}}
\def\L{{\mathcal L}}
\def\cK{{\mathcal K}}
\def\cJ{{\mathcal J}}
\def\cS{{\mathcal S}}
\def\bH{{\bf H}}
\def\bP{{\bf P}}
\def\bQ{{\bf Q}}
\def\bE{{\bf E}}
\def\bT{{\bf T}}
\def\W{W}
\def\Be{L_\infty}
\def\cR{{\mathcal R}}
\def\eps{\varepsilon}
\def\3{{\ss}}
\def\slim{s-\lim_}
\def\capa{{\mathrm{Cap}}}
\def\supp{{\mathrm{supp}}}
\def\esssup{{\mathrm{ess\,sup}}}
\def\absconv{{\mathrm{absconv}}}
\def\dom{{\mathrm{dom}}}
\def\loc{\mathrm{loc}}
\def\hs{half-space}
\def\HIC{$\HH^\infty$-calculus}
\def\BIP{{\mathrm{BIP}}}                                                 
\def\BUC{{\mathrm{BUC}}}
\def\BC{{\mathrm{BC}}}
\def\MR{{\mathcal{MR}}}
\def\const{{\mathrm{const\,}}}
\def\Re{{\mathrm{Re}}}
\def\re{{\mathrm{Re}}}
\def\Im{{\mathrm{Im}}}
\def\im{{\mathrm{Im}}}
\def\dd{{\mathrm d}}
\def\e{{\mathrm{e}}}
\def\id{{\mathrm{id}}}
\def\sb{{\mathrm{sb}}}
\def\FM{{\mathrm{FM}}}
\def\ME{{\mathrm{M}}}
\def\hperp{{^{_\perp}}}
\def\HIC{{$H^\infty$-calculus}}
\def\hW{\widehat{W}}
\def\hu{\hat{u}}
\def\hv{\hat{v}}
\def\hw{\hat{w}}
\def\hsigma{\hat{\sigma}}
\def\hf{\hat{f}}
\def\hh{\hat{h}}
\def\hg{\hat{g}}
\def\dR{\dot{\R}}
\def\tu{\tilde{u}}
\def\tc{\tilde{c}}
\def\tp{\tilde{p}}
\def\tf{\tilde{f}}
\def\th{\tilde{h}}
\def\tg{\tilde{g}}
\def\tv{\tilde{v}}
\def\ta{\tilde{a}}
\def\ty{\widetilde{y}}
\def\bv{\bar{v}}
\def\bw{\bar{w}}
\def\tsigma{\tilde{\sigma}}
\def\hphi{\hat{\phi}}
\newcommand{\essinf}[1]{{\mathrm{ess}}\!\inf_{\!\!\!\!\!\!\!\!\! #1}}
\newcommand{\fn}{\footnote}
\def\mdt{mixed derivative theorem}
\def\div{{\mathrm {div\,}}}
\def\bsigma{\bar{\sigma}}
\def\brho{\bar{\rho}}
\def\rcv{\W^1_p(J;L^p(\R^{n+1}_+))
        \cap L_p(J;\W^2_p(\R^{n+1}_+))}
\def\rcvtp{\W^1_p(J;L^p(\dR^{n+1}))
        \cap L_p(J;\W^2_p(\dR^{n+1}))}
\def\rcs{W^{3/2-1/2p}_p(J;L_p(\R^n))
        \cap \W^1_p(J;W^{1-1/p}_p(\R^n))
        \cap L_p(J;W^{2-1/p}_p(\R^n))}
\def\rcf{L_p(J;L_p(\R^{n+1}_+))}
\def\rcftp{L_p(J;L_p(\dR^{n+1}))}
\def\rch{W^{1/2-1/2p}_p(J;L_p(\R^n))
        \cap L_p(J;W^{1-1/p}_p(\R^n))}
\def\rcvi{W^{2-2/p}_p(\R^{n+1}_+)}
\def\rcvitp{W^{2-2/p}_p(\dR^{n+1})}
\def\rcsi{W^{2-2/p}_p(\R^n)}
\newcommand{\ab}{&\hskip-2mm}

\def\en{{\talloblong}}
\def\hookd{\stackrel{_d}{\hookrightarrow}}
\def\hook{{\hookrightarrow}}
\def\vt{{\vartheta}}
\def\ovt{{\overline{\vartheta}}}
\def\THE{tornado-hurricane equations}
\def\THO{tornado-hurricane operator}
\def\NSE{Navier-Stokes equations}
\def\SO{Stokes operator}
\def\HHP{Helmholtz projection}
\def\HHD{Helmholtz decomposition}
\def\HOL{\mathrm{HOL}}
\def\la{{\langle}}
\def\ra{{\rangle}}
\def\vphi{{\varphi}}
\def\vdp{{\{k:\ \alpha^v_k=0\}}}
\def\vndp{{\{k:\ \alpha^v_k\neq0\}}}
\def\tdp{{\{k:\ \alpha^\vt_k=0\}}}
\def\tndp{{\{k:\ \alpha^\vt_k\neq0\}}}
\def\sD{{\mathscr D}}
\def\hsD{{\widehat{\mathscr D}}}
\def\sL{{\mathscr L}}
\def\sR{{\mathscr R}}
\def\sT{{\mathscr T}}
\def\sLis{{\mathscr L}_{is}}
\def\PPr{{\PP_{\!\rho}}}
\def\ou{{\overline{u}}}
\def\oq{{\overline{q}}}
\def\oU{{\overline{U}}}
\def\th{{T\!H}}
\def\ttau{{\tilde \tau}}
\def\hH{{\widehat{H}}}
\def\cD{{\mathcal D}}
\def\vp{{\varphi}}
\def\Hic{{\mathcal{H}^\infty}}

\hyphenation{Lipschitz}

%\newcommand\fn[1]{}

%**************************************************************************
\sloppy
%**************************************************************************
\title[Stable and unstable flow regimes for active fluids]
{Stable and unstable flow regimes for active fluids in the periodic setting}

\author[C. Bui]{Christiane Bui}
\address{
         Heinrich-Heine-Uni\-ver\-sit\"at D\"usseldorf\\
	Mathematisches Institut, Angewandte Analysis\\
         40204 D\"usseldorf, Germany\\
	 email:\ christiane.bui@hhu.de}
%\email{Christiane.Bui@uni-duesseldorf.de}

\author[C. Gesse]{Christian Gesse}
\address{
         Heinrich-Heine-Uni\-ver\-sit\"at D\"usseldorf\\
	Mathematisches Institut, Angewandte Analysis\\
         40204 D\"usseldorf, Germany\\
	 email:\ christian.gesse@hhu.de}
%\email{Christian.Gesse@uni-duesseldorf.de}

\author[J. Saal]{J\"urgen Saal}
\address{
         Heinrich-Heine-Uni\-ver\-sit\"at D\"usseldorf\\
	Mathematisches Institut, Angewandte Analysis\\
         40204 D\"usseldorf, Germany\\
	 %phone:\ +49 211 81-11366\\
	 email:\ juergen.saal@hhu.de}
%\email{juergen.saal@hhu.de}

%\date{\today}
%\thanks{class-stefan75.tex}
%\thispagestyle{empty}
%\pagestyle{myheadings}
%\markboth{\today}{\today}
%\setlength{\parindent}{0mm}
%\setlength{\parskip}{10cm}
\parskip0.5ex plus 0.5ex minus 0.5ex
%\bibliographystyle{plain}
%%\bibliographystyle{alpha}
%\bibliographystyle{amsalpha}
%\setcounter{page}{3}
%\listoffigures

\begin{abstract}

Depending on the involved physiobiological parameters, stable or
unstable behavior in active fluids is observed. In this paper a 
rigorous analytical justification of (in-)stability within the
corresponding regimes is given. 
In particular, occuring instability for the manifold of ordered polar states 
caused by self-propulsion is proved. This represents the prerequisite 
for active turbulence patterns as observed in a number of applications. 
The approach is carried out in the periodic
setting and is based on the generalized principle of linearized
(in)-stability related to normally stable and normally hyperbolic equilibria. 
\end{abstract}
\maketitle
%%%%%%%%%%%%%%%%%%%%%%%%%%%%%%%%%%%%

{\bf Keywords.} Living fluids, active turbulence, generalized Navier-Stokes
equations, periodic setting, well-posedness, stability, 
%{\bf 2000 Mathematics Subject Classification.}
%Primary ; Secondary 
 
%\tableofcontents

%\doublespacing

%%%%%%%%%%%%%%%%%%%%%%%%%%%%%%%%%%%%%%%%%%%%%%%%%%%%%%%%%%%%%%%%%%%%%%%%%%%%%
\section{Introduction}
%%%%%%%%%%%%%%%%%%%%%%%%%%%%%%%%%%%%%%%%%%%%%%%%%%%%%%%%%%%%%%%%%%%%%%%%%%%%%
%%%%%%%%%%%%%%%%%%%%%%%%%%%%%%%%%%%%%%%%%%%%%%%%%%%%%%%%%%%%%%%%%%%%%%%%%%%%%

A minimal hydrodynamic model to describe the bacterial velocity in
the case of highly concentrated bacterial suspensions with negligible density
fluctuations is given as
\begin{equation}
	\label{eqn:min-hyd-mod}
	\begin{array}{rl}
		v_t+\lambda_0v\cdot\nabla v & =  f-\nabla
		p+\lambda_1\nabla|v|^2-(\alpha+\beta|v|^2)v+\Gamma_0\Delta
		v-\Gamma_2 \Delta^2v,\\
		\mbox{div}\,v & =  0,\\
		v(0)&=v_0,
	\end{array}
\end{equation}
see \cite{Wensink-et-al:Meso-scale-turbulence}.
Here $v$ is the bacterial velocity field and
$p$ the (scalar) pressure and $\lambda_0, \lambda_1, \alpha, \beta, \Gamma_0$ and $\Gamma_2$ are real parameters.

In \cite{zls2016} a first rigoros analytical approach to (\ref{eqn:min-hyd-mod})
in $L^2(\R^n)$ is presented. There, depending on the values of the 
involved parameters,
results on (in-) stability of the disordered isotropic steady state
and ordered polar steady
states, see \eqref{disorderedstate} and \eqref{orderedstate}, are derived. 
A formal stability analysis based on the standard wave ansatz coming
to the same conclusions is already performed in 
\cite{Wensink-et-al:Meso-scale-turbulence}. This formal analysis,
however, cannot be rigorously confirmed by the approach in $L^2(\R^n)$
as given in \cite{zls2016}, just by the fact that the wave
ansatz is not an $L^2(\R^n)$-function. Therefore in \cite{bls2019} 
an approach to (\ref{eqn:min-hyd-mod}) in spaces of Fourier transformed 
Radon measures $\mathrm{FM}(\R^n)$ is developed. It gives the same
outcome on (in)-stability as in \cite{zls2016}. Moreover, it confirms
the formal stability analysis in 
\cite{Wensink-et-al:Meso-scale-turbulence}, since the space 
$\mathrm{FM}(\R^n)$ contains wave functions such as 
$\exp(ik\cdot x+\sigma t)$ required for the wave ansatz.

A drawback of the approaches presented in \cite{zls2016} and \cite{bls2019}
concerns the instability of the ordered polar states. The results in 
\cite{zls2016} and \cite{bls2019} prove that a single polar state
for specific regimes of the involved parameters is unstable. The polar
states, however, form the manifold $B_{\alpha,\beta}$, i.e., the sphere with
radius $\sqrt{-\alpha/\beta}$ centered at the origin, 
see \eqref{orderedstate}. The results
in \cite{zls2016} and \cite{bls2019} do not clarify the question, if 
a solution still could converge to the manifold $B_{\alpha,\beta}$, even
though each single polar state on $B_{\alpha,\beta}$ displays unstable
behavior. On the other hand, instability of the manifold of
polar states is the prerequisite for active turbulence, and thus
especially interesting with regard to applications, see
\cite{Wensink-et-al:Meso-scale-turbulence,Thampi,Doostmohammadi,Dogic}.

The purpose of this note is to develop an approach in $L^2$ in the
periodic setting. This allows for 
applying the concepts of normally stable and 
normally hyperbolic equilibria as, e.g.,
provided by \cite[Theorem 2.1 and Theorem~6.1]{psz2009} 
and \cite[Theorem~5.3.1 and Theorem~5.5.1]{moving_interfaces}.
Utilizing the latter
result we can precisely characterize instability of the manifold 
$B_{\alpha,\beta}$. In fact, compactness yields
discrete spectrum of the involved linear operators. A crucial
point then is to prove that zero is a semi-simple eigenvalue of the
linearization at each equilibrium on the manifold, which will turn out
to be true for the polar states. It should be noticed 
that this strategy is not
applicable in $L^2(\R^n)$ and $\mathrm{FM}(\R^n)$, essentially by the 
fact that the spectra are continuous in those settings. 
Comparing the three approaches, it appears that the approach in 
the periodic setting performed here seems
the most suitable one. This is also underlined by the fact that
the generalized Navier-Stokes system (\ref{eqn:min-hyd-mod})
was augmented even further by
including further higher derivatives in the velocities entering into
the stress tensor, see \cite{Slomka_EPJST_2015}. In presence of 
a boundary, this required additional higher order boundary conditions.
From the physical point of view it seems not at all 
clear what are suitable boundary conditions to be imposed.
Remaining in the periodic setting it is no problem to include also
higer order terms in the velocity. Indeed, we anticipate that the 
analysis performed here also applies to this more general case, provided
the highest order term has the correct sign.

Note that the model (\ref{eqn:min-hyd-mod})
was originally proposed 
in \cite{Wensink-et-al:Meso-scale-turbulence} 
and then further considered in
\cite{Dunkel_Heidenreich_PRL_2013,Dunkel_NJP_2013}.
For $\lambda_0=1$, $\lambda_1 = \alpha=\beta = \Gamma_2 = 0$ 
and $\Gamma_0 >0$, it reduces to the classical 
incompressible Navier-Stokes equations 
in $n$ spatial dimensions. 
For non-vanishing 
$\lambda_1,\alpha,\beta,\Gamma_2$ system
(\ref{eqn:min-hyd-mod}) by now is one of the standard models 
to describe active turbulence
at low Reynolds number \cite{Oza}. It can also be derived from more 
microscopic descriptions \cite{Klapp} and it was quantitatively confirmed in
suspensions of living biological systems \cite{Wensink-et-al:Meso-scale-turbulence,Kaiser1,Dogic,Beppu} and synthetic microswimmers 
\cite{Sagues}. Note that active turbulence was also suggested as a power source for various microfluidic applications 
\cite{Kaiser1,Kaiser2,Kaiser3,Thampi}. 
We refer to those papers and to \cite{zls2016} for
a more detailed description of the physics behind the additional occuring
terms. 

We organized this note as follows: Section~\ref{sec_perspace} 
collects basic facts on periodic spaces in the $L^2$-setting.
In Section~\ref{subsec_wp} we analize the linearized system about 
the relevant equilibria. Here precise statements on the spectra 
and the corresponding asymptotic behavior of the semigroups are
established. Furthermore, we prove global-in-time well-posedness for
system (\ref{eqn:min-hyd-mod}) in the strong setting.
Section~\ref{subsec_stab} then concerns the nonlinear active turbulence.
In Subsection~\ref{subsec_nt_ds} we first consider the stability behavior
of the disordered state depending on the involved parameters and 
relying on the linear stability analysis. Subsection~\ref{subsec_nt_psns}
deals with the stable regime for the manifold of the polar states. It will be
proved to be normally stable in this case. The most important result 
related to active turbulence 
then is given in Subsection~\ref{subsec_nt_psnh}. There the manifold of 
ordered polar states is proved to be normally hyperbolic in the
unstable regime of the parameters.

%%%%%%%%%%%%%%%%%%%%%%%%%%%%%%%%%%%%%%%%%%%%%%%%%%%%%%%%%%%%%%%%%%%%%%%%%%%%%
\section{Periodic Sobolev spaces}\label{sec_perspace}
%%%%%%%%%%%%%%%%%%%%%%%%%%%%%%%%%%%%%%%%%%%%%%%%%%%%%%%%%%%%%%%%%%%%%%%%%%%%%
%%%%%%%%%%%%%%%%%%%%%%%%%%%%%%%%%%%%%%%%%%%%%%%%%%%%%%%%%%%%%%%%%%%%%%%%%%%%%
We start with some basic notation. Let $\Omega \subseteq \R^n$ be a domain.
By $L^p(\Omega,X)$ for $1 \leq p \leq \infty$ 
we denote the standard Bochner-Lebesgue space with values
in a Banach space $X$. As usual we equip $L^p(\Omega,X)$ with the standard
Lebesgue space norm
	\begin{align*}
	\| u \|_{L^p(\Omega,X)} =
	\left ( \int_\Omega \|u(x)\|_X^p dx \right )^{1/p}
	\end{align*}
for $1 \leq p < \infty$ with the usual modification if $p = \infty$. 
For $k \in \N$ and $1 \leq p \leq \infty$ the
Sobolev space $W^{k,p}(\Omega,X)$ of $k$ times weakly differentiable 
functions is equipped with the norm
	\begin{align*}
	\|u\|_{W^{k,p}(\Omega,X)} =
		\left ( \sum_{|\alpha| \leq k} 
		\|\partial^\alpha u\|_{L^p(\Omega,X)}^p \right )^{1/p}.
	\end{align*}
If $p = 2$ then we write $H^k(\Omega,X) := W^{k,2}(\Omega,X)$. 
Next, we define periodic Sobolev spaces. 
Let $L > 0$ be arbitrary but fixed from now on. We set $Q_n := [0,L]^n$. 
The function space which corresponds to 
periodic boundary conditions is $L^2_\pi(Q_n, \R^n)$. It is defined 
as the completion of $C^\infty_\pi(Q_n)$ 
with respect to the $L^2(Q_n,\R^n)$ norm, where
    \begin{align*}
    C^k_\pi(Q_n) &:=
    \left \{ f \in C^{k}(Q_n,\R^n) : \partial^\alpha f|_{x_j = 0}
    = \partial^\alpha f|_{x_j = L} \: 
    		\forall |\alpha| \leq k \right \},\\
    C^\infty_\pi(Q_n) &:= \bigcap_{k = 0}^\infty C^k_\pi(Q_n).
    \end{align*}
Here, $C^k(\Omega,X)$ denotes the space of $k$ times continuously
differentiable functions with values in a Banach space $X$.
To simplify the notation we set 
$L^2(Q_n) := L^2_\pi(Q_n) := L^2_\pi(Q_n, \R^n)$. 
By \cite[Proposition 3.2.1]{grafakos} it follows that the definitions 
of $L^2_\pi(Q_n,\R^n)$ and $L^2(Q_n,\R^n)$ as a standard Lebesgue space
are equivalent.

One advantage of working in $L^2(Q_n)$ is the fact that we can employ 
Fourier series. For a detailed introduction to Fourier series 
and further important properties of $L^2(Q_n)$ we refer to 
\cite[Chapter 3]{grafakos} and 
\cite[Chapter 5.10]{infinite_dimensional_dynamical_systems}.
Let $f \in L^2(Q_n)$. The Fourier coefficient $\hat{f}(m)$ for 
$m=(m_1,...,m_n) \in \Z^n$ is defined as the integral
    \[
    \hat{f}(m):= \mathcal{F}f(m) :=
    \frac{1}{L^n} \int_{Q_n} f(x) e^{-2\pi i m x / L} \: dx
    \]
where $mx = \sum_{k=1}^n m_k x_k$ denotes the scalar product in 
$\R^n$. By using integration by parts one can verify the identity
    \begin{align}
    \label{per_fourierseries_deriv}
    \widehat{\partial^\alpha f}(m) 
    = \left (\frac{2\pi i}{L} \right )^{|\alpha|}
    m^\alpha \hat{f}(m)
    \end{align}
if $f$ is smooth enough, i.e., $f \in C_\pi^{|\alpha|}(Q_n)$,
$m \in \Z^n$ and $\alpha \in \N_0^n$. Let $f, g  \in L^2(Q_n)$. 
Then the scalar product in $L^2(Q_n)$ is defined as
    \begin{align*}
    (f,g)_{2,\pi} := \frac{1}{L^n} \int_{Q_n} f(x) \overline{g(x)} \, dx
    \end{align*}
where the subscript $\pi$ indicates that we are in the periodic $L^2$ space.
Some important properties of the Fourier series and 
$L^2(Q_n)$ functions are listed below (\cite[Proposition 3.2.7]{grafakos}).
%%%%%%%%%%%%%%%%%%%%%%%%%%%%%%%%%%%%%%%%%%%%%%%%%%%%%%%%%%%%%%%%%%%%%%%%%%%%%
    \begin{proposition}
    \label{per_important_prop}
    Let $f, g \in L^2(Q_n)$ be arbitrary. The following properties hold in $L^2(Q_n)$:
\begin{enumerate}
        \item Plancherel theorem:
            \begin{align*}
            \| f \|_{L^2(Q_n)}^2 = \sum_{m \in \Z^n} |\hat{f}(m)|^2.
            \end{align*}
        \item Parseval's identity:
            \begin{align*}
            (f,g)_{2,\pi} 
            = \frac{1}{L^n} \int_{Q_n} f(x) \overline{g(x)} \, dx
            = \sum_{m \in \Z^n} \hat{f}(m) \overline{\hat{g}(m)}.
            \end{align*}
        \item The function $f$ 
            can be represented as the $L^2(Q_n)$-limit of trigonometric
            polynomials, i.e.,
            \begin{align*}
            f = \sum_{m \in \Z^n} \hat{f}(m) e^{2 \pi i m \cdot / L}.
            \end{align*}
	    	\end{enumerate}
    \end{proposition}
%%%%%%%%%%%%%%%%%%%%%%%%%%%%%%%%%%%%%%%%%%%%%%%%%%%%%%%%%%%%%%%%%%%%%%%%%%%%%
We will consider (\ref{eqn:min-hyd-mod}) with periodic boundary conditions 
in $L^2(Q_n)$.
For this purpose we define relevant periodic Sobolev spaces for $k \in \N$:
    \begin{align*}
    H^k_\pi(Q_n) 
    &:= \left \{
         u = \sum_{m \in \Z^n} \hat{u}(m) e^{2\pi i m \cdot/L} :
        \hat{u}(m) = \overline{\hat{u}(-m)},
        \|u\|_{\tilde{H}_\pi^k(Q_n)} < \infty
        \right \}\\
    &= \left \{ u \in H^k(Q_n) : \partial^\alpha u|_{x_j = 0} 
    		= \partial^\alpha u|_{x_j = L} \:
    			(|\alpha| < k, j = 1,...,n) \right \}\\
    &= \overline{C^\infty_\pi(Q_n)}^{H^k(Q_n)}
    \end{align*}
where the norm above is defined as
    \begin{align*}
    \|u\|_{\tilde{H}_\pi^k(Q_n)}^2 
    := \sum_{m \in \Z^n} \left | \left (1+\left ( 
    		\frac{2\pi}{L} \right )^{k}|m|^{k} \right ) \hat{u}(m) \right|^2,
    \end{align*}
see \cite[Chapter 5.10]{infinite_dimensional_dynamical_systems}.
If $u \in H^k_\pi(Q_n)$ and $\alpha \in \N_0^n$ with $|\alpha| \leq k$, 
then the derivative $\partial^\alpha u$ can be written as the 
$L^2(Q_n)$-limit
    \begin{align*}
    \partial^\alpha u = \sum_{k \in \Z^n} \widehat{\partial^\alpha u}(k)
        e^{2\pi i k \cdot / L}
        = \sum_{k \in \Z^n} \left ( \frac{2\pi i}{L}\right )^{|\alpha|} 
        		k^\alpha \hat{u}(k)
        e^{2\pi i k \cdot / L},
    \end{align*}
where we used the fact that the identity in (\ref{per_fourierseries_deriv}) 
also holds for $u \in H^k_\pi(Q_n)$. It obviously follows that the 
$\| \cdot \|_{\tilde{H}_\pi^k(Q_n)}$ and the $\| \cdot \|_{H_\pi^k(Q_n)}$ 
norms are equivalent, where 
    \begin{align*}
    \|u\|_{H_\pi^k(Q_n)}^2
    := \sum_{|\alpha|\leq k} \sum_{m \in \Z^n} \left | \left ( 
    		\frac{2 \pi}{L}\right )^{|\alpha|} m^\alpha \hat{u}(m)\right |^2,
    \end{align*}
by the Plancherel theorem.
Periodic Sobolev spaces of fractional powers are defined in the
canonical way:
For $s \geq 0$ we set
    \begin{align*}
    H^s_\pi(Q_n) 
    &= \left \{
         u = \sum_{m \in \Z^n} \hat{u}(m) e^{2\pi i m \cdot/L} :
        \hat{u}(m) = \overline{\hat{u}(-m)},
        \|u\|_{H_\pi^s(Q_n)} < \infty
        \right \}
    \end{align*}
where
    \begin{align*}
    \|u\|_{H_\pi^s(Q_n)}^2 
    := \sum_{m \in \Z^n} \left (1+\left( \frac{2\pi}{L}\right)^2 |m|^2 \right )^{s/2} |\hat{u}(m)|^2,
    \end{align*}
and it is straightforward to see that for $s \in \N$ the two definitions
for Sobolev spaces coincide.

Finally, let $m : \Z^n \to \C^{n \times n}$ be a function. We define 
$T_m : D(T_m) \subseteq L^2(Q_n) \to L^2(Q_n)$ as the $L^2(Q_n)$-limit
    \begin{align*}
    T_mf := \sum_{k \in \Z^n} m(k) \hat{f}(k) e^{2\pi i k \cdot / L}
    \end{align*}
for a function $f \in D(T_m)$ where
	\begin{align*}
	D(T_m) := \left \{ f \in L^2(Q_n) : \|T_m f\|_{L^2(Q_n)}^2 = \sum_{k \in \Z^n} |m(k)\hat{f}(k)|^2 < \infty \right \}.
	\end{align*}
Note that $T_m$ is well-defined and a bounded operator if $m$ is a bounded 
function by the Plancherel theorem. Then $m$ is called a Fourier multiplier 
on $L^2(Q_n)$.

%%%%%%%%%%%%%%%%%%%%%%%%%%%%%%%%%%%%%%%%%%%%%%%%%%%%%%%%%%%%%%%%%%%%%%%%%%%%%
\section{Linear stability and well-posedness}\label{subsec_wp}
%%%%%%%%%%%%%%%%%%%%%%%%%%%%%%%%%%%%%%%%%%%%%%%%%%%%%%%%%%%%%%%%%%%%%%%%%%%%%
%%%%%%%%%%%%%%%%%%%%%%%%%%%%%%%%%%%%%%%%%%%%%%%%%%%%%%%%%%%%%%%%%%%%%%%%%%%%%
In the following we will consider two different, physically 
relevant stationary solutions
of (\ref{eqn:min-hyd-mod}):
	\begin{equation}\label{disorderedstate}
		(v,p)=(0,p_0),
	\end{equation}
which corresponds to a \textit{disordered isotropic state}, and, 
if $\alpha<0$, the set of equilibra corresponds to the manifold 
of \textit{globally ordered polar states}:
	\begin{equation}\label{orderedstate}
		(v,p)=(V,p_0),
	\end{equation}
where $V\in B_{\alpha,\beta}
:=\{x\in \mathbb R^n:\ |x|=\sqrt{-\alpha/\beta}\}$, i.e., $V$ 
denotes a constant vector with arbitrary orientation
and fixed swimming speed $|V|=\sqrt{-\alpha/\beta}$. In both cases the 
pressure $p_0$ is a constant.

In order to cover all situations corresponding to the above steady states, 
as in \cite{zls2016}, we consider the following generalization of 
(\ref{eqn:min-hyd-mod}):
	\begin{align}
		\label{eqn:min-hyd-mod-trans}
		\begin{array}{rcl}
			u_t+\lambda_0\left[(u+V)\cdot\nabla\right] u 
			+(M+\beta|u|^2)u \qquad \\[0.25em] -\Gamma_0\Delta
			u+\Gamma_2\Delta^2u+\nabla q& = &
			  f+N(u), \\[0.25em]
		\mbox{div}\,u & =&  0,\\[0.25em]    
		u(0)&=&u_0,
		\end{array}
	\end{align}
where $q=p-\lambda_1 |v|^2$, $M\in\R^{n\times n}$ is a symmetric matrix and
$N(u)=\sum_{j,k}a_{jk}u^ju^k$ with $(a_{jk})_{j,k=1}^n \in \R^{n \times n}$ 
defines a nonlinearity of second order. 
Regarding the occuring parameters we assume
	\begin{equation}\label{domparameters}
		\lambda_0,\lambda_1,\Gamma_0,\alpha\in\R,
		\qquad \Gamma_2,\beta>0,
	\end{equation}
throughout this paper. From (\ref{eqn:min-hyd-mod-trans}) we obtain
the equation corresponding to the disordered state (\ref{disorderedstate}) by setting
	\begin{equation}\label{valuesds}
		V=0,\quad M=\alpha I,\quad N(u)=0,
	\end{equation}
for $u=v$ where $I$ denotes the identity matrix and $\alpha$ is a scalar. By setting
	\begin{equation}\label{valuesos}
		V\in B_{\alpha,\beta},\quad M=2\beta VV^T,\quad 
		N(u)=-\beta|u|^2V-2\beta(u\cdot V)u
	\end{equation}
we obtain the system for $u=v-V$ corresponding to the ordered polar state (\ref{orderedstate}).
Furthermore, space dimension is always assumed to be $n=2$ or $n=3$.
%%%%%%%%%%%%%%%%%%%%%%%%%%%%%%%%%%%%%%%%%%%%%%%%%%%%%%%%%%%%%%%%%%%%%%%%%%%%%
\subsection{The linearized system}\label{sec_lin}
%%%%%%%%%%%%%%%%%%%%%%%%%%%%%%%%%%%%%%%%%%%%%%%%%%%%%%%%%%%%%%%%%%%%%%%%%%%%%
%%%%%%%%%%%%%%%%%%%%%%%%%%%%%%%%%%%%%%%%%%%%%%%%%%%%%%%%%%%%%%%%%%%%%%%%%%%%%
In this subsection we consider the linearized system
    \begin{equation}
    \label{lflin}
    \begin{array}{r@{\ =\ }ll}
        u_t + \lambda_0 (V \cdot \nabla) u + Mu 
        - \Gamma_0 \Delta u + \Gamma_2 \Delta^2 u + \nabla q
        & f &\text{in } (0, \infty) \times Q_n,\\
        \mbox{div}\,u & 0 &\text{in } (0, \infty) \times Q_n,\\
        u(0) & u_0 &\text{in } Q_n
    \end{array}
    \end{equation}
with periodic boundary conditions
	\begin{align*}
	\partial^\alpha u|_{x_j = 0} 
	= \partial^\alpha u|_{x_j = L}
	\quad \text{for } |\alpha| < 4, j = 1,..., n.
	\end{align*}
First we introduce the Helmholtz-Weyl projection on $L^2(Q_n)$.
Setting the Fourier multiplier $\sigma_P : \Z^n \to \C^{n \times n}$ as
$\sigma_P(m) = I - mm^T/|m|^2$ for $m \neq 0$ and $\sigma_P(0) = I$, the 
Helmholtz-Weyl projection $P : L^2(Q_n) \to L^2_\sigma(Q_n)$ is defined as
    \begin{align*}
    Pu = \sum_{m \in \Z^n} \sigma_P(m) \hat{u}(m) e^{2 \pi i m \cdot / L}
    \qquad (u \in L^2(Q_n)).
    \end{align*}
The projection $P$ induces the Helmholtz decomposition 
\[
	L^2(Q_n) = L^2_\sigma(Q_n) \oplus G_2(Q_n),
\]
where
    \begin{align*}
    L^2_\sigma(Q_n) &:=
    \left \{ u \in L^2(Q_n) : 
        \hat{u}(m) = \overline{\hat{u}(-m)}, m \cdot \hat{u}(m) = 0
        \: \forall \: m \in \Z^n
    \right \}, \\
        G_2(Q_n) &:=
    \left \{
        u = \nabla g\in L^2(Q_n) : g \in L^1_{loc}(Q_n)
    \right \}.
    \end{align*}
Note that $P$ is also a projection on $H^{k}_\pi(Q_n)$ and that 
$P(H^{k}_\pi(Q_n))=H^k_\pi(Q_n) \cap L^2_\sigma(Q_n)$. From this
we obtain the following interpolation result.
%%%%%%%%%%%%%%%%%%%%%%%%%%%%%%%%%%%%%%%%%%%%%%%%%%%%%%%%%%%%%%%%%%%%%%%%%%%%%
    \begin{lemma}
    \label{interpolation_periodic_spaces}
    Let $\theta \in [0,1]$ and $k \in \N$. Then we have
        \begin{align*}
        [L^2_\sigma(Q_n),H^k_\pi(Q_n) \cap L^2_\sigma(Q_n)]_\theta
        = H^{\theta k}_\pi(Q_n) \cap L^2_\sigma(Q_n),
        \end{align*}
    where $[\cdot,\cdot]_\theta$ denotes the complex interpolation
    functor, see \cite{triebel}.
    \end{lemma}
%%%%%%%%%%%%%%%%%%%%%%%%%%%%%%%%%%%%%%%%%%%%%%%%%%%%%%%%%%%%%%%%%%%%%%%%%%%%%
   \begin{proof}
   The assertion follows from the fact that $P$ is a projection
   on the interpolated spaces and
   \cite[Theorem 1.2.4]{triebel}.
   \end{proof}
%%%%%%%%%%%%%%%%%%%%%%%%%%%%%%%%%%%%%%%%%%%%%%%%%%%%%%%%%%%%%%%%%%%%%%%%%%%%%
We define the operator associated to (\ref{lflin}) as
    \begin{align*}
    %\label{per_def_op}
    \begin{split}
    A_{LF}u &:= \lambda_0 (V \cdot \nabla) u 
        + PMu - \Gamma_0 \Delta u + \Gamma_2 \Delta^2u,\\
    D(A_{LF}) &:= H^4_\pi(Q_n) \cap L^2_\sigma(Q_n),
    \end{split}
    \end{align*}
and the corresponding Fourier symbol as
    \begin{align}
    \label{per_op_symbol}
    \sigma_{A_{LF}}(\ell)
    := \Gamma_2 \left ( \frac{2\pi}{L} \right )^4 |\ell|^4
        + \Gamma_0 \left ( \frac{2\pi}{L} \right )^2 |\ell|^2
        + \lambda_0 \left ( \frac{2\pi i}{L} \right ) (V \cdot \ell)
        + \sigma_P(\ell)M
    \end{align}
for $\ell \in \Z^n$. Note that thanks to $\Gamma_2 > 0$ we immediately see that the operator
    \begin{align}
    \label{per_def_ah}
    A_{SH} u := \Gamma_2 \Delta^2 u, \quad
    D(A_{SH}) := H^4_\pi(Q_n) \cap L^2_\sigma(Q_n),
    \end{align}
is selfadjoint. Hence $\omega + A_{SH}$ is selfadjoint and positive 
for $\omega > 0$.
As a consequence $\omega + A_{SH}$ admits 
a bounded $H^\infty$-calculus on $L^2_\sigma(Q_n)$
with $H^\infty$-angle $\phi_{\omega + A_{SH}}^\infty = 0$.
See e.g.\ \cite{haase2006} for an introduction to the notion of a bounded
$H^\infty$-calculus. Applying perturbation theorems
we obtain the following result:
%%%%%%%%%%%%%%%%%%%%%%%%%%%%%%%%%%%%%%%%%%%%%%%%%%%%%%%%%%%%%%%%%%%%%%%%%%%%%
	\begin{proposition}
    \label{per_hinfty}
    There exists an $\omega > 0$ such that $\omega + A_{LF}$
    admits a bounded
    $H^\infty$-calculus on $L^2_\sigma(Q_n)$ with $H^\infty$-angle $\phi_{\omega + A_{LF}}^\infty = 0$.
    \end{proposition}
%%%%%%%%%%%%%%%%%%%%%%%%%%%%%%%%%%%%%%%%%%%%%%%%%%%%%%%%%%%%%%%%%%%%%%%%%%%%%
    	\begin{proof}
    	This follows from the fact that
    		\begin{align*}
	    %\label{per_def_b}
	    \begin{split}
	    Bu &:= \lambda_0 ( V \cdot \nabla) u + PMu - \Gamma_0 \Delta u
	    \end{split}
	    \end{align*}
	is a perturbation of lower order. The assertion then follows from 
	\cite[Proposition 13.1]{Kunstmann2004}.
    	\end{proof}
%%%%%%%%%%%%%%%%%%%%%%%%%%%%%%%%%%%%%%%%%%%%%%%%%%%%%%%%%%%%%%%%%%%%%%%%%%%%%
As a consequence $A_{LF}$ enjoys maximal $L^p$-regularity on intervals
$(0,T)$ with $T < \infty$ and $-A_{LF}$ is the generator of an analytic
$C_0$-semigroup on $L^2_\sigma(Q_n)$:
%%%%%%%%%%%%%%%%%%%%%%%%%%%%%%%%%%%%%%%%%%%%%%%%%%%%%%%%%%%%%%%%%%%%%%%%%%%%%
	\begin{proposition}
	\label{per_op_maxreg}
	Let $T \in (0,\infty)$. For $f \in L^2((0,T),L^2_\sigma(Q_n))$
	and $u_0 \in H^2_\pi(Q_n) \cap L^2_\sigma(Q_n)
	=\bigl(L^2_\sigma(Q_n),H^4_\pi(Q_n)\cap
	L^2_\sigma(Q_n)\bigr)_{1/2,2}$ 
	there exists a unique
	solution $(u,q)$ of (\ref{lflin}) such that
		\begin{align*}
		&\|u\|_{H^1((0,T),L^2_\sigma(Q_n))}
		+ \|u\|_{L^2((0,T),H^4_\pi(Q_n))}
		+ \|\nabla q\|_{L^2((0,T),L^2_\sigma(Q_n))}\\
		&\leq C(T) \left ( \|f\|_{L^2((0,T),L^2_\sigma(Q_n))}
		+ \|u_0\|_{H^2_\pi(Q_n)}\right ),
		\end{align*}
	where $C(T) > 0$ is independent of $u,q,u_0,f$.
	\end{proposition}
%%%%%%%%%%%%%%%%%%%%%%%%%%%%%%%%%%%%%%%%%%%%%%%%%%%%%%%%%%%%%%%%%%%%%%%%%%%%%
 It is straightforward to verify the identity
	\begin{align}
	\label{semigroup_rep}
	\exp(-tA_{LF})v = \sum_{\ell \in \Z^n}
	\exp\bigl(-t\sigma_{A_{LF}}(\ell)\bigr)
	\hat{v}(\ell)e^{2\pi i \ell \cdot /L}
	\end{align}
for $v \in L^2_\sigma(Q_n)$. Using this representation we can characterize
linear (in)stability for both stationary states (\ref{disorderedstate}) and
(\ref{orderedstate}). First we examine the disordered isotropic state
(\ref{disorderedstate}). In this case we set $A_d := A_{LF}$ where
$V = 0$ and $M = \alpha I$. Then $P$ commutes with $M$ and $PMu = \alpha u$
for $u \in D(A_d)$. The corresponding Fourier symbol is given as
    \begin{align*}
    %\label{per_op_symbol_disordered}
    \sigma_{A_d}(\ell) := \Gamma_2 \left ( \frac{2\pi}{L} \right )^4 |\ell|^4
        + \Gamma_0 \left ( \frac{2\pi}{L} \right )^2 |\ell|^2
        + \alpha
    \qquad (\ell \in \Z^n).
    \end{align*}
Using representation (\ref{semigroup_rep}) we immediately obtain 
%%%%%%%%%%%%%%%%%%%%%%%%%%%%%%%%%%%%%%%%%%%%%%%%%%%%%%%%%%%%%%%%%%%%%%%%%%%%%
    \begin{proposition}
    \label{per_semigroup_stab}
    Let $\Gamma_2 > 0$ and $\Gamma_0, \alpha \in \R$. Then the semigroup
    $(\exp(-tA_d))_{t \geq 0}$ corresponding to the disordered state
    (\ref{disorderedstate}) is
        \begin{enumerate}
        \item stable if $\sigma_{A_d} \geq 0$;
        \item exponentially stable if $\sigma_{A_d} \geq \delta > 0$;
        \item exponentially unstable if there exists some $\ell_0 \in \Z^n$
        such that $\sigma_{A_d}(\ell_0) < 0$.
        \end{enumerate}
    \end{proposition}
%%%%%%%%%%%%%%%%%%%%%%%%%%%%%%%%%%%%%%%%%%%%%%%%%%%%%%%%%%%%%%%%%%%%%%%%%%%%%
%    \begin{proof}
%    At first we consider the case $\sigma_{A_d} \geq \delta > 0$. Using Plancherel's
%    theorem (Theorem~\ref{per_important_prop}(1)) we see that
%        \begin{align*}
%        \| \exp(-tA_d) v \|_{L^2(Q_n)}^2
%        \leq \sum_{\ell \in \Z^n} e^{-2t\delta} |\hat{v}(\ell)|^2
%        = e^{-2t\delta} \|v\|^2_{L^2(Q_n)}
%        \end{align*}
%    holds for $v \in L^2_\sigma(Q_n)$    
%    which yields exponential stability of the semigroup $(\exp(-tA_d))_{t \geq 0}$. On the
%    other hand, if there exists some $\ell_0 \in \Z^n$ such that
%    $\sigma_{A_d}(\ell_0) < 0$ we obtain
%        \begin{align*}
%        \| \exp(-tA_d) v \|_{L^2(Q_n)}^2
%        \geq e^{-2t\sigma_{A_d}(\ell_0)} |\hat{v}(\ell_0)|^2 \to \infty
%        \end{align*}
%    as $t \to \infty$. In this case the semigroup $(\exp(-tA_d))_{t \geq 0}$ is exponentially unstable.
%    To see stability we note that
%        \begin{align*}
%        \| \exp(-tA_d) v \|_{L^2(Q_n)}^2
%        = \sum_{\ell \in \Z^n} e^{-2t\sigma_{A_d}(\ell)} |\hat{v}(\ell)|^2
%        \leq \| v\|_{L^2(Q_n)}^2
%        \end{align*}
%    since $\sigma_{A_d} \geq 0$.
%    \end{proof}
%%%%%%%%%%%%%%%%%%%%%%%%%%%%%%%%%%%%%%%%%%%%%%%%%%%%%%%%%%%%%%%%%%%%%%%%%%%%%
Now, we can characterize linear stability in terms of the involved parameters by considering the behaviour of the symbol $\sigma_{A_d}$ for a certain choice of $\Gamma_0$ and $\alpha$. For this purpose we substitute $z = |\ell|^2 \geq 0$:
    \begin{align*}
    p(z) := \Gamma_2 \left ( \frac{2\pi}{L} \right )^4 z^2
        + \Gamma_0 \left ( \frac{2\pi}{L} \right )^2 z
        + \alpha
    \end{align*}
and we see that $p$ describes a parabola. The intersection points are given as
    \begin{align*}
    z_{\pm} = \frac{-\Gamma_0}{\Gamma_2 \left ( \frac{2\pi}{L} \right)^2}
        \left ( \frac{1}{2} \pm \sqrt{\frac{1}{4} - \frac{\alpha \Gamma_2}{\Gamma_0^2}} \right ),
    \end{align*}
if $\Gamma_0 \neq 0$ and as
    \begin{align*}
    z_\pm = \pm \sqrt{\frac{-\alpha}{\Gamma_2} \left (\frac{L}{2\pi} \right ) ^4}
    \end{align*}
if $\Gamma_0 = 0$.
If $\Gamma_0 < 0$ and $4\alpha > \Gamma_0^2/\Gamma_2$ or $\Gamma_0 \geq 0$ and $\alpha > 0$ then there exists some $\delta > 0$ such that $\sigma_{A_d} \geq \delta$. 
Using these observations and Proposition~\ref{per_semigroup_stab} we
obtain the following concreter classification of stability:
%%%%%%%%%%%%%%%%%%%%%%%%%%%%%%%%%%%%%%%%%%%%%%%%%%%%%%%%%%%%%%%%%%%%%%%%%%%%%
    \begin{corollary}
    \label{per_linstab_disordered}
    Let $\Gamma_2 > 0$. If $\Gamma_0 < 0$ and $4\alpha > \Gamma_0^2/\Gamma_2$
    or if $\Gamma_0 \geq 0$ and $\alpha > 0$ then the semigroup 
    $(\exp(-tA_d))_{t \geq 0}$ generated
    by $-A_d$ is exponentially stable. 
    To be precise, the semigroup $(\exp(-tA_d))_{t \geq 0}$
    corresponding to the disordered state (\ref{disorderedstate}) is
        \begin{enumerate}
         \item stable, if $\Gamma_0 < 0$ and $4\alpha \geq
            \Gamma_0^2/\Gamma_2$ or if $\Gamma_0 \geq 0$ and $\alpha \geq 0$;
        \item exponentially stable, if $\Gamma_0 < 0$ and $4\alpha > \Gamma_0^2/\Gamma_2$ or if $\Gamma_0 \geq 0$ and $\alpha > 0$
        or if $\Gamma_0 < 0$ and $4\alpha = \Gamma_0^2/\Gamma_2$ with
        $|\ell|^2 \neq -\frac{\Gamma_0}{2\Gamma_2} \left(\frac{L}{2\pi}\right)^2$ for all $\ell \in \Z^n$.
        \end{enumerate}
    \end{corollary}
%%%%%%%%%%%%%%%%%%%%%%%%%%%%%%%%%%%%%%%%%%%%%%%%%%%%%%%%%%%%%%%%%%%%%%%%%%%%%
    \begin{proof}
    From the observations above and by Proposition~\ref{per_semigroup_stab}, 
    exponential stability for the cases $\Gamma_0 < 0$ and
    $4\alpha > \Gamma_0^2/\Gamma_2$ and $\Gamma_0 \geq 0$ and $\alpha > 0$
    follows immediately. 
    If $\Gamma_0 \geq 0$ and $\alpha = 0$ then
    $\sigma_{A_d}(0) = 0$ and $\sigma_{A_d}(\ell) \geq 0$ otherwise such
    that we can apply Proposition~\ref{per_semigroup_stab} to see stability.
    If $\Gamma_0 < 0$ and $4\alpha = \Gamma_0^2/\Gamma_2$ then the
    symbol $\sigma_{A_d}$ simplifies as
        \begin{align*}
        \sigma_{A_d}(\ell) 
            &= \left ( 
            \sqrt{\Gamma_2} \left ( \frac{2\pi}{L} \right )^2 |\ell|^2
            + \frac{\Gamma_0}{2\sqrt{\Gamma_2}}
            \right )^2
        \end{align*}
    for $\ell \in \Z^n$. We have $\sigma_{A_d} > 0$ if and only if
        \begin{align*}
        |\ell|^2 \neq - \frac{\Gamma_0}{2 \Gamma_2}
            \left ( \frac{L}{2\pi} \right )^2
        \end{align*}
    for all $\ell \in \Z^n$. Again the assertion follows from
    Proposition~\ref{per_semigroup_stab}.
    \end{proof}
%%%%%%%%%%%%%%%%%%%%%%%%%%%%%%%%%%%%%%%%%%%%%%%%%%%%%%%%%%%%%%%%%%%%%%%%%%%%%
Next, we consider the ordered polar state (\ref{orderedstate}). Here, we
set $A_o := A_{LF}$ where $V \in B_{\alpha,\beta}$ and $M = 2\beta VV^T$.
The Fourier symbol is then given as
    \begin{align*}
    \sigma_{A_o}(\ell)
        := \Gamma_2 \left ( \frac{2\pi}{L} \right )^4 |\ell|^4
        &+ \Gamma_0 \left ( \frac{2\pi}{L} \right )^2 |\ell|^2
        + \lambda_0 \left ( \frac{2\pi i}{L} \right ) (V \cdot \ell)\\
        &+ 2\beta\sigma_P(\ell)VV^T\sigma_P(\ell).
    \end{align*}
If $\Gamma_0 \geq 0$ then
	\begin{align*}
	\Re \, \sigma_{A_o}(\ell) =
	\Gamma_2 \left ( \frac{2\pi}{L} \right )^4 |\ell|^4
        + \Gamma_0 \left ( \frac{2\pi}{L} \right )^2 |\ell|^2
        + 2\beta\sigma_P(\ell)VV^T\sigma_P(\ell)
	 \in \R^{n \times n}
	\end{align*}
is positive semi-definite
for every $\ell \in \Z^n$ and even positive definite
for $\ell \neq 0$. This follows from the fact that
$\sigma_P(\ell)VV^T\sigma_P(\ell)$ is positive semi-definite.
Then we can estimate the norm of the semigroup as 
    \begin{align*}
    \|\exp(-tA_o)v\|_{L^2(Q_n)}^2
    \leq |\hat{v}(0)|^2 + \sum_{\ell \in \Z^n \backslash \{0\}} |e^{-t\sigma_{A_o}(\ell)}|^2 |\hat{v}(\ell)|^2
    \end{align*}
to see stability.
Conversely, if $\Gamma_0 < 0$ and if there exists $0 \neq \ell_0 \in \Z^n$ such that
    \begin{align}
    \label{per_linstab_ordered_eq}
    	\Gamma_2 \left ( \frac{2\pi}{L} \right )^2 |\ell_0|^2
    	+ \Gamma_0
        < 0
    \end{align}
then for $n = 3$ we can choose  $x \in \R^n \setminus \{ 0 \}$ with $x \perp V$, $x \perp \ell_0$ such that $x^T \Re \ \sigma_{A_o}(\ell_0)x < 0$. In case of $n = 2$, due to
\begin{align*}
	- \frac{2 \beta |V\cdot \hat{v}(\ell)|^2}{|\hat{v}(\ell)|^2} \in [2\alpha, 0]
\end{align*} 
(see also (\ref{nonlinear-hyperbolicity-condition}) later)
we assume the existence of $0 \neq \ell_0 \in \Z^n$ such that
\begin{align}\label{per_linstab_ordered_instab_cond}
	\Gamma_2 \left ( \frac{2\pi}{L} \right )^2 |\ell_0|^2
	+ \Gamma_0 < 2 \alpha.
\end{align}
This implies $x^T \Re \ \sigma_{A_o}(\ell_0)x < 0$ for 
some $x \in \R^n \setminus \{ 0 \}$ such that $x\perp \ell_0$.
Consequently, the matrix $\Re \ \sigma_{A_o}(\ell_0) \in \R^{n \times n}$ 
is negative semi-definite or indefinite. Hence the growth bound of 
$(\exp(-tA_o))_{t\ge 0}$ is strictly positive.
We obtain the following result on stability:
%%%%%%%%%%%%%%%%%%%%%%%%%%%%%%%%%%%%%%%%%%%%%%%%%%%%%%%%%%%%%%%%%%%%%%%%%%%%%
    \begin{proposition}
    \label{per_linstab_ordered}
    Let $\Gamma_2 > 0$. Then the semigroup $(\exp(-tA_o))_{t \geq 0}$ corresponding to the ordered
    polar state (\ref{orderedstate}) is
		\begin{enumerate}
		\item stable if $\Gamma_0 \geq 0$;
		\item exponentially unstable if $\Gamma_0 < 0$ and 
		\begin{enumerate}
			\item[(i)] if there exists some $0 \neq \ell_0 \in \Z^n$ such that (\ref{per_linstab_ordered_instab_cond}) holds for $n=2$;
			\item[(ii)] if there exists some $0 \neq \ell_0 \in \Z^n$ such that (\ref{per_linstab_ordered_eq}) holds for $n=3$.
		\end{enumerate}
		\end{enumerate}	    
    \end{proposition}
%%%%%%%%%%%%%%%%%%%%%%%%%%%%%%%%%%%%%%%%%%%%%%%%%%%%%%%%%%%%%%%%%%%%%%%%%%%%%

%%%%%%%%%%%%%%%%%%%%%%%%%%%%%%%%%%%%%%%%%%%%%%%%%%%%%%%%%%%%%%%%%%%%%%%%%%%%%
\begin{remark}
(a)\ It is worthwhile to compare at this point the situation to the
continuous case considered in \cite[Section~3.1]{zls2016} and
\cite[Section~3.1]{bls2019}. There a $\xi\in\R^n$
such that $\xi$ is parallel to $V$ can always be found, 
which in general is not possible in
the discrete case. Consequently, also for dimension $n=2$ a nontrivial 
$x\in\R^2$ satisfying $x\perp V$ and $x\perp \xi$ and giving instability 
does exist. By this fact, in \cite{zls2016} and \cite{bls2019} 
the more restrictive condition \eqref{per_linstab_ordered_instab_cond}
for $n=2$ does not appear.\\[1mm]
(b)\ Note that for $n=2$ condition
(\ref{per_linstab_ordered_instab_cond}) imposes no restrictions
regarding the analysis of nonlinear instability considered here, cf. condition (\ref{cond_norm_hyper}) in Theorem \ref{per_nonlin_hyperbolic_ordered}.
\end{remark}
%%%%%%%%%%%%%%%%%%%%%%%%%%%%%%%%%%%%%%%%%%%%%%%%%%%%%%%%%%%%%%%%%%%%%%%%%%%%%

%%%%%%%%%%%%%%%%%%%%%%%%%%%%%%%%%%%%%%%%%%%%%%%%%%%%%%%%%%%%%%%%%%%%%%%%%%%%%
%%%%%%%%%%%%%%%%%%%%%%%%%%%%%%%%%%%%%%%%%%%%%%%%%%%%%%%%%%%%%%%%%%%%%%%%%%%%%
\subsection{Global well-posedness}\label{subsec_nlwp}
%%%%%%%%%%%%%%%%%%%%%%%%%%%%%%%%%%%%%%%%%%%%%%%%%%%%%%%%%%%%%%%%%%%%%%%%%%%%%
%%%%%%%%%%%%%%%%%%%%%%%%%%%%%%%%%%%%%%%%%%%%%%%%%%%%%%%%%%%%%%%%%%%%%%%%%%%%%
In this section we quote the result on global well-posedness of
(\ref{eqn:min-hyd-mod-trans}). Here the proof is omitted, since it 
is completely analogous to the proof in \cite[Section 3.2]{zls2016}.
%%%%%%%%%%%%%%%%%%%%%%%%%%%%%%%%%%%%%%%%%%%%%%%%%%%%%%%%%%%%%%%%%%%%%%%%%%%%%
    \begin{theorem}[global well-posedness]
    \label{per_gwp}
    Let $\Gamma_2, \beta > 0$ and $\Gamma_0, \alpha, \lambda_0 \in \R$
    and $T \in (0, \infty)$. Let the initial value $u_0 \in H^2_\pi(Q_n)
    \cap L^2_\sigma(Q_n)$ and an exterior force 
    $f \in L^2((0,T),L^2_\sigma(Q_n))$ be given. 
    Then there exists a unique pair $(u,q)$
    with
        \begin{align*}
        u &\in H^1((0,T), L^2_\sigma(Q_n)) \cap L^2((0,T),H^4_\pi(Q_n)),\\
        \nabla q &\in L^2((0,T), L^2(Q_n)),
        \end{align*}
    solving (\ref{eqn:min-hyd-mod-trans}) for periodic boundary conditions.
    \end{theorem}
%%%%%%%%%%%%%%%%%%%%%%%%%%%%%%%%%%%%%%%%%%%%%%%%%%%%%%%%%%%%%%%%%%%%%%%%%%%%%
	\begin{remark}
	Note that, in contrast to the classical Navier-Stokes equations,
	the convective term in (\ref{eqn:min-hyd-mod-trans}) is
	dominated by the linear fourth order term. Consequently,
	standard energy techniques lead to global strong
	solvability, see \cite[Section 3.2]{zls2016} for the details.
\end{remark}
%%%%%%%%%%%%%%%%%%%%%%%%%%%%%%%%%%%%%%%%%%%%%%%%%%%%%%%%%%%%%%%%%%%%%%%%%%%%%

%%%%%%%%%%%%%%%%%%%%%%%%%%%%%%%%%%%%%%%%%%%%%%%%%%%%%%%%%%%%%%%%%%%%%%%%%%%%%
\section{Nonlinear stability and turbulence}\label{subsec_stab}
%%%%%%%%%%%%%%%%%%%%%%%%%%%%%%%%%%%%%%%%%%%%%%%%%%%%%%%%%%%%%%%%%%%%%%%%%%%%%
%%%%%%%%%%%%%%%%%%%%%%%%%%%%%%%%%%%%%%%%%%%%%%%%%%%%%%%%%%%%%%%%%%%%%%%%%%%%%
In this section we study nonlinear stability of the stationary states of (\ref{eqn:min-hyd-mod}). 
We will consider the disordered state (\ref{disorderedstate}) and 
the manifold of ordered polar states (\ref{orderedstate}) separately. To see stability 
we will apply the generalized principle of linearized stability, cf.\
\cite[Theorem 2.1]{psz2009} 
and \cite[Theorem 5.3.1]{moving_interfaces}, and energy methods; to see 
instability we will use the principle of normally hyperbolic equilibria
\cite[Theorem 6.1]{psz2009} resp.\ 
\cite[Theorem 5.5.1]{moving_interfaces} and Henry's instability theorem
\cite[Corollary 5.1.6]{henry}.

%%%%%%%%%%%%%%%%%%%%%%%%%%%%%%%%%%%%%%%%%%%%%%%%%%%%%%%%%%%%%%%%%%%%%%%%%%%%%
\subsection{The disordered state}\label{subsec_nt_ds}
%%%%%%%%%%%%%%%%%%%%%%%%%%%%%%%%%%%%%%%%%%%%%%%%%%%%%%%%%%%%%%%%%%%%%%%%%%%%%
%%%%%%%%%%%%%%%%%%%%%%%%%%%%%%%%%%%%%%%%%%%%%%%%%%%%%%%%%%%%%%%%%%%%%%%%%%%%%

We start with the following auxiliary result.
%%%%%%%%%%%%%%%%%%%%%%%%%%%%%%%%%%%%%%%%%%%%%%%%%%%%%%%%%%%%%%%%%%%%%%%%%%%%%
    \begin{lemma}
    \label{per_nonlin_henry_assump}
    Let
    $H(u) := \beta P|u|^2 u + \lambda_0 P(u \cdot \nabla)u - PN(u)$.
    Then we have
    $H \in C^1(H^\eta_\pi(Q_n) \cap L^2_\sigma(Q_n), L^2_\sigma(Q_n))$ for
    $\eta \geq 5/4$ and $H$ can be estimated as
        \begin{align*}
        \|H(u)\|_{L^2(Q_n)} \leq C \|u\|_{H^\eta_\pi(Q_n)}^2
        \qquad (\|u\|_{H^\eta_\pi(Q_n)} \leq 1).
        \end{align*}
    \end{lemma}
%%%%%%%%%%%%%%%%%%%%%%%%%%%%%%%%%%%%%%%%%%%%%%%%%%%%%%%%%%%%%%%%%%%%%%%%%%%%%
	\begin{proof}
	The proof of the lemma follows verbatim the lines of the proof
	of \cite[Lemma 4]{zls2016}.
	\end{proof}
%%%%%%%%%%%%%%%%%%%%%%%%%%%%%%%%%%%%%%%%%%%%%%%%%%%%%%%%%%%%%%%%%%%%%%%%%%%%%
Now, we consider the disordered isotropic state
(\ref{disorderedstate}). Suppose $(u,q)$ is the global solution to
(\ref{eqn:min-hyd-mod-trans}) from Theorem~\ref{per_gwp}. Then we have
    \begin{align*}
    %\label{per_nonlinstab_dis_eq}
    u_t + \Gamma_2 \Delta^2 u - \Gamma_0 \Delta u 
    + \lambda_0 (u \cdot \nabla) u + (\alpha + \beta |u|^2)u
    + \nabla q = 0.
    \end{align*}
Testing the equation above with $u$ we obtain
    \begin{align}
    \label{per_nonlinstab_dis_eq_2}
    \begin{split}
    \frac{1}{2} \frac{d}{dt} \|u(t)\|_{L^2(Q_n)}^2
    &+ \Gamma_2 \|\Delta u(t)\|_{L^2(Q_n)}^2
    + \Gamma_0 \|\nabla u(t)\|_{L^2(Q_n)}^2\\
    &+ \alpha \|u(t)\|_{L^2(Q_n)}^2
    + \beta \|u(t)\|_{L^4(Q_n)}^4\\
    &= 0
    \end{split}
    \end{align}
for $t > 0$.
%%%%%%%%%%%%%%%%%%%%%%%%%%%%%%%%%%%%%%%%%%%%%%%%%%%%%%%%%%%%%%%%%%%%%%%%%%%%%
    \begin{theorem}
    \label{per_nonlin_stab_disordered}
    Let $\Gamma_2, \beta > 0$ and $\Gamma_0, \alpha, \lambda_0 \in \R$.
    Then the
    disordered state (\ref{disorderedstate}) is nonlinearly
        \begin{enumerate}
        \item stable in $L^2_\sigma(Q_n)$ if $\Gamma_0 \geq 0$ and $\alpha \geq 0$ or if $\Gamma_0 < 0$
        and $4\alpha \geq \Gamma_0^2/\Gamma_2$;
        \item (globally) exponentially stable in $L^2_\sigma(Q_n)$ if $\Gamma_0 \geq 0$ and
        $\alpha > 0$ or if $\Gamma_0 < 0$ and $4\alpha > \Gamma_0^2/\Gamma_2$
        or if $\Gamma_0 < 0$ and $4\alpha = \Gamma_0^2/\Gamma_2$ and
        $|\ell|^2 \neq \frac{-\Gamma_0}{2\Gamma_2} \left ( \frac{L}{2\pi}\right )^2$ for all $\ell \in \Z^n$;
        \item unstable in $H^\gamma_\pi(Q_n) \cap L^2_\sigma(Q_n)$ for $\gamma \in [5/16,1)$ if there
        exists some $\ell_0 \in \Z^n$ such that $\sigma_{A_d}(\ell_0) < 0$.
        \end{enumerate}
    \end{theorem}
%%%%%%%%%%%%%%%%%%%%%%%%%%%%%%%%%%%%%%%%%%%%%%%%%%%%%%%%%%%%%%%%%%%%%%%%%%%%%
    \begin{proof}
    If $\Gamma_0 \geq 0$ and $\alpha \geq 0$
    we immediately obtain from (\ref{per_nonlinstab_dis_eq_2})
        \begin{align*}
        \frac{1}{2} \frac{d}{dt} \|u(t)\|_{L^2(Q_n)}^2
        + \alpha \|u(t)\|_{L^2(Q_n)}^2 \leq 0,
        \end{align*}
    since $\Gamma_2, \Gamma_0, \beta \geq 0$. Applying Gronwall's lemma
    we deduce
        \begin{align*}
        \|u(t)\|_{L^2(Q_n)}^2 \leq e^{-2\alpha t} \|u_0\|_{L^2(Q_n)}^2,
        \end{align*}
    which shows that the disordered state (\ref{disorderedstate}) is exponentially stable
    if $\alpha > 0$ and stable if $\alpha = 0$. 
    To see exponential
    stability in the other case ($\Gamma_0 < 0$ and $4\alpha > \Gamma_0^2/\Gamma_2$) we drop the $\beta$ term 
    in (\ref{per_nonlinstab_dis_eq_2}) and apply Plancherel's theorem
    to obtain
        \begin{align*}
        \frac{1}{2} \frac{d}{dt} \|u(t)\|_{L^2(Q_n)}^2
        &+ \sum_{\ell \in \Z^n} \left ( 
        \Gamma_2 \left ( \frac{2\pi}{L} \right )^4 |\ell|^4
        + \Gamma_0 \left ( \frac{2\pi}{L}\right )^2 |\ell|^2
        + \alpha \right ) |\widehat{u(t)}(\ell)|^2 \\
        &\leq 0.
        \end{align*}
    We note that we can estimate the Fourier symbol in the series as
        \begin{align*}
        %\label{per_nonlinstab_disordered_est_1} 
        \Gamma_2 \left ( \frac{2\pi}{L} \right )^4 |\ell|^4
        + \Gamma_0 \left ( \frac{2\pi}{L}\right )^2 |\ell|^2 
        + \alpha > \delta
        \end{align*}
    for some $\delta > 0$ by applying Young's inequality with
    $\varepsilon^2 = 2\Gamma_2 / |\Gamma_0|$:
        \begin{align*}
        \Gamma_2 \left ( \frac{2\pi}{L} \right )^4 |\ell|^4
        &+ \Gamma_0 \left ( \frac{2\pi}{L}\right )^2 |\ell|^2 
        + \alpha
        \geq \alpha - \frac{\Gamma_0^2}{4\Gamma_2} > 0.
        \end{align*}
    Hence in this case we set 
    $\delta := \alpha - \Gamma_0^2/4\Gamma_2 > 0$
    to see
        \begin{align*}
        \frac{1}{2} \frac{d}{dt} \|u(t)\|_{L^2(Q_n)}^2
        + \delta \|u(t)\|_{L^2(Q_n)}^2
        \leq 0.
        \end{align*}
    Again the application of Gronwall's lemma  
    yields exponential stability. Next, 
    if $\Gamma_0 < 0$ and $4\alpha = \Gamma_0^2/\Gamma_2$
    then $\sigma_{A_d}$
    simplifies as
        \begin{align*}
        \sigma_{A_d}(\ell)
        = \left ( \sqrt{\Gamma_2} \left ( \frac{2\pi}{L} \right )^2 |\ell|^2 + \frac{\Gamma_0}{2\sqrt{\Gamma_2}}\right )^2
        \end{align*}
    for $\ell \in \Z^n$ and we have $\sigma_{A_d} > 0$
    if and only if 
    $|\ell|^2 \neq \frac{-\Gamma_0}{2\Gamma_2} \left ( \frac{L}{2\pi}\right )^2$ for every $\ell \in \Z^n$. 
    Then there exists $\delta > 0$ such
    that $\sigma_{A_d} > \delta$
    and using the same arguments as before we infer exponential
    stability. 
    Conversely, if there exists some $\ell_0 \in \Z^n$ such that 
    $|\ell_0|^2 = \frac{-\Gamma_0}{2\Gamma_2} \left ( \frac{L}{2\pi}\right )^2$
    then $\sigma_{A_d} \geq 0$ and the disordered state (\ref{disorderedstate})
    remains stable.
    
    Using \cite[Corollary 5.1.6]{henry}
    we can show instability: In the notation of \cite[Corollary 5.1.6]{henry}
    we have $x_0 = 0$, $x = u$, $A = A_d$ and $f(u) = g(u) = H(u)$. Then $-A_d$
    generates an analytic semigroup and there
    exists an $\omega > 0$ such that $\omega + A_d$ admits a
    bounded $H^\infty$-calculus,
    see Proposition~\ref{per_hinfty}. Then we have
        \begin{align*}
        D(A_d^\gamma) = [L^2_\sigma(Q_n), H^4_\pi(Q_n) \cap L^2_\sigma(Q_n)]_\gamma
        = H^{4\gamma}_\pi(Q_n) \cap L^2_\sigma(Q_n)
        \end{align*}
    for $\gamma \in [0,1]$ which follows from Lemma~\ref{interpolation_periodic_spaces}. 
    Furthermore, under the assumptions on $\sigma_{A_d}$ we know that 
    the disordered state is exponentially
    unstable for the linear system such that
        \begin{align*}
        \sigma(-A_d) \cap \{ z \in \C : \Re \, z > 0 \} \neq \emptyset,
        \end{align*}
    see Proposition~\ref{per_semigroup_stab}. 
    With Lemma~\ref{per_nonlin_henry_assump} all conditions of 
    \cite[Corollary 5.1.6]{henry} hold for $\gamma \in [5/16,1)$ and the 
    disordered state (\ref{disorderedstate}) is unstable in this case.
    \end{proof}
%%%%%%%%%%%%%%%%%%%%%%%%%%%%%%%%%%%%%%%%%%%%%%%%%%%%%%%%%%%%%%%%%%%%%%%%%%%%%

%%%%%%%%%%%%%%%%%%%%%%%%%%%%%%%%%%%%%%%%%%%%%%%%%%%%%%%%%%%%%%%%%%%%%%%%%%%%%
\subsection{Ordered polar states: normal stability}\label{subsec_nt_psns}
%%%%%%%%%%%%%%%%%%%%%%%%%%%%%%%%%%%%%%%%%%%%%%%%%%%%%%%%%%%%%%%%%%%%%%%%%%%%%
%%%%%%%%%%%%%%%%%%%%%%%%%%%%%%%%%%%%%%%%%%%%%%%%%%%%%%%%%%%%%%%%%%%%%%%%%%%%%

Now, we consider the manifold $B_{\alpha,\beta}$ of ordered polar 
states (\ref{orderedstate}) for the stable regime.
Let $V\in B_{\alpha,\beta}$ and $A_o$ be the corresponding linear
operator as defined in Section~\ref{sec_lin}.

An equilibrium $V \in B_{\alpha,\beta}$ is called normally stable,
cf.\ \cite[Theorem 5.3.1]{moving_interfaces}, if
	\begin{enumerate}
	\item[(i)] near $V$ the set of equilibria $B_{\alpha,\beta}$ is a $C^1$-manifold in $H^4_\pi(Q_n) \cap L^2_\sigma(Q_n)$ of dimension $m \in \N$;
	\item[(ii)] the tangent space $T_V B_{\alpha,\beta}$ for $B_{\alpha,\beta}$ at $V$ is isomorphic to $N(A_o)$;
	\item[(iii)] $0$ is a semisimple eigenvalue of $A_o$, i.e., $L^2_\sigma(Q_n) = N(A_o) \oplus R(A_o)$;
	\item[(iv)] $\sigma(A_o) \backslash \{0\} \subseteq \{z \in \C : \Re \ z > 0\}$.
	\end{enumerate}
The equilibrium $V$ is called normally hyperbolic, 
cf.\ \cite[Theorem 5.5.1]{moving_interfaces}, 
if the conditions (i) - (iii) hold and
	\begin{enumerate}
	\item[(iv)'] $\sigma(A_o) \cap i \R = \{0\}$ and $\sigma_u := \sigma(A_o) \cap \C_- \neq \emptyset$.
	\end{enumerate}
For the stable regime we will show exponential stability by applying 
the generalized principle of 
linearized stability \cite[Theorem 2.1]{psz2009}, 
\cite[Theorem 5.3.1]{moving_interfaces}.
For this purpose, we first prove
%%%%%%%%%%%%%%%%%%%%%%%%%%%%%%%%%%%%%%%%%%%%%%%%%%%%%%%%%%%%%%%%%%%%%%%%%%%%%
    \begin{lemma}
    \label{per_nonlinstab_prop_a0}
    Let
    $\Gamma_0 \geq 0$. Then $0$
    is a semisimple eigenvalue of $A_o$, i.e., $N(A_o) \oplus R(A_o) = L^2_\sigma(Q_n)$. Furthermore, we can characterize the spectrum of $A_o$ in the
    following way: The spectrum of $A_o$ only consists of eigenvalues
    and is discrete. Additionally,
        \begin{align*}
        \sigma(A_o) \subseteq \{ \lambda \in \C : \emph{Re} \, \lambda > 0 \} \cup
            \{ 0 \}.
        \end{align*}
    \end{lemma}
%%%%%%%%%%%%%%%%%%%%%%%%%%%%%%%%%%%%%%%%%%%%%%%%%%%%%%%%%%%%%%%%%%%%%%%%%%%%%
    \begin{proof}
    First we note that for $\lambda \in \rho(A_o)$
    \begin{align*}
    (\lambda - A_o)^{-1} : 
    L^2_\sigma(Q_n) \to D(A_o) 
    = H^4_\pi(Q_n) \cap L^2_\sigma(Q_n) \overset{c}{\hookrightarrow} L^2_\sigma(Q_n)
    \end{align*}
is compact by the Rellich-Kondrachov theorem (\cite[Theorem A.4, Corollary A.5]{infinite_dimensional_dynamical_systems}).
	Then $\sigma(A_o)$ is discrete and $\sigma(A_o) = \sigma_p(A_o)$ where
	$\sigma_p(A_o)$ denotes the point spectrum of $A_o$.

    Next, we show that $0$ is a semisimple eigenvalue of $A_o$. 
    Let
    $u \in H^4_\pi(Q_n) \cap L^2_\sigma(Q_n)$ be a constant such that
    $u$ is perpendicular to $V$. Then
        \begin{align*}
        A_o u
        = \Gamma_2 \Delta^2 u - \Gamma_0 \Delta u
            + \lambda_0 (V \cdot \nabla) u
            + 2 \beta P V V^T u
        = 0
        \end{align*} 
    and $0$ is
    an eigenvalue of $A_o$. We want to characterize $N(A_o)$ more
    precisely. From the argumentation above we immediately have
        \begin{align*}
        \{ u \in H^4_\pi(Q_n) \cap L^2(Q_n) : u \: 
        \text{constant and} \: u \perp V \} \subseteq N(A_o).
        \end{align*}
    To prove the converse inclusion we take $u \in N(A_o)$. Then $A_o u = 0$
    and testing the equation with $u$ we obtain
        \begin{align*}
        0
        = (\Gamma_2 \Delta^2 u,u)_{2,\pi} - (\Gamma_0 \Delta u, u)_{2,\pi}
            + (\lambda_0 (V \cdot \nabla)u, u)_{2,\pi}
            + 2 \beta (P V V^T u, u)_{2,\pi}.
        \end{align*}
    Taking the real part and applying integration by parts we derive
        \begin{align*}
        0 = \Gamma_2 \| \Delta u \|_{L^2(Q_n)}^2 
            + \Gamma_ 0 \|\nabla u\|_{L^2(Q_n)}^2
            + 2 \beta \|V \cdot u\|_{L^2(Q_n)}^2,
        \end{align*}
    since the $\lambda_0$ term is skew-symmetric.
    By assumptions $\Gamma_2, \beta > 0$ and $\Gamma_0 \geq 0$ we conclude
        \begin{align*}
        \|\Delta u\|_{L^2(Q_n)}^2
        = \|V \cdot u\|_{L^2(Q_n)}^2 = 0,
        \end{align*}
    hence $u$ is constant and perpendicular to $V$ since
    	\begin{align*}
    	\|\Delta u\|_{L^2(Q_n)}^2 = \sum_{\ell \in \Z^n} |\ell|^2 |\hat{u}(\ell)|^2 = 0.
    	\end{align*}
    Consequently,
        \begin{align*}
        N(A_o) =
        \{ u \in H^4_\pi(Q_n) \cap L^2_\sigma(Q_n) : u \: \text{constant and} \:
         u \perp V \}.
        \end{align*}
        
    Next, we show the decomposition $N(A_o) \oplus R(A_o) = L^2_\sigma(Q_n)$.
    We define the following map
        \begin{align*}
        S : L^2_\sigma(Q_n) \to L^2_\sigma(Q_n),
        \qquad Su := \frac{1}{L^n} \int_{Q_n} S_* u(x) dx,
        \end{align*}
    where $S_* : L^2_\sigma(Q_n) \to L^2_\sigma(Q_n)$ is the map given by $S_*u(x) = (I - VV^T/|V|^2)u(x)$,
    where $I$ denotes the identity matrix in $n$ dimensions. 
    First we note that if $u \in L^2_\sigma(Q_n)$ then $Su$ is constant
    and $Su \in L^2_\sigma(Q_n)$.
    It is straightforward to prove that
    $S$ is a projection such that there exists a decomposition $S(L^2_\sigma(Q_n))
    \oplus (I-S)(L^2_\sigma(Q_n)) = L^2_\sigma(Q_n)$. Then we need to show
    that $S(L^2(Q_n)) = N(A_o)$ and $(I-S)(L^2_\sigma(Q_n)) = R(A_o)$.
    
    First we claim $N(A_o) = S(L^2_\sigma(Q_n))$. 
    To see the inclusion
    $S(L^2_\sigma(Q_n)) \subseteq N(A_o)$ we assume $u \in S(L^2_\sigma(Q_n))$.
    Then $u = Su$  is constant as already mentioned and we observe
        \begin{align*}
        V^T u &= V^T Su
        &= \frac{1}{L^n} \int_{Q_n} V^T u(x) dx
            - \frac{1}{L^n} \int_{Q_n} \frac{1}{|V|^2} V^T V V^T u(x) dx
        = 0,
        \end{align*}
    hence $u=Su$ is perpendicular to $V$ which yields $u \in N(A_o)$. To see the converse
    inclusion we take $u \in N(A_o)$. Then $u$ is constant and perpendicular to
    $V$. We obtain
        \begin{align*}
        Su 
        &= \frac{1}{L^n} \int_{Q_n} u \, dx
            - \frac{1}{L^n} \int_{Q_n} \frac{1}{|V|^2} VV^T u \, dx\\
        &= u \left ( \frac{1}{L^n} \int_{Q_n} dx \right ) = u,
        \end{align*}
    hence $u \in S(L^2_\sigma(Q_n))$ and the claim is proved.
    Since $L^2_\sigma(Q_n)$ is a Hilbert space and $S$ is a selfadjoint projection
    it is well known that $L^2_\sigma(Q_n) = S(L^2_\sigma(Q_n))
    \oplus (I-S)(L^2_\sigma(Q_n))$ is an orthogonal decomposition. If
    we take $u \in D(A_o)$ and show that $A_o u$ is perpendicular to any $w
    \in N(A_o)$ then $R(A_o) \subseteq (I-S)(L^2_\sigma(Q_n))$:
        \begin{align*}
        (A_o u, w)_{2,\pi}
        &= \Gamma_2 (\Delta u, \Delta w)_{2,\pi} 
        		+ \Gamma_0 (\nabla u, \nabla w)_{2,\pi}\\
            &- \lambda_0 (u, (V \cdot \nabla)w)_{2,\pi}
            + 2 \beta (V^Tu, V^T w)_{2,\pi}\\
        &= 0,
        \end{align*}
    because $w$ is constant and perpendicular to $V$. 
    Since $A_o$ has compact resolvent it follows from \cite[Corollary 1.19]{en00} that the spectral value $0$
    is a pole of the resolvent. Then by \cite[Remark A.2.4]{lunardi}
    it suffices to
    show that
        \begin{align*}
        N(A_o) = N(A_o^2)
        \end{align*}
    to prove that $0$ is a semisimple eigenvalue of $A_o$. The inclusion $N(A_o)
    \subseteq N(A_o^2)$ is obvious. To see the
    converse inclusion we take $u \in N(A_o^2)$. Then $A_o^2 u = 0$ such that
    $A_o u \in N(A_o)\cap R(A_o)=\{0\}$ by what we just proved.
    Thus, $N(A_o^2) = N(A_o)$.
    
    The last assertion $\sigma(A_o) \subseteq \{\lambda \in \C : \Re \, \lambda > 0 \} \cup \{0\}$
    follows from Proposition~\ref{per_linstab_ordered}(i) and the fact that $-A_o$ generates a bounded holomorphic $C_0$-semigroup in this case.
    \end{proof}
%%%%%%%%%%%%%%%%%%%%%%%%%%%%%%%%%%%%%%%%%%%%%%%%%%%%%%%%%%%%%%%%%%%%%%%%%%%%%
Now, we are in position to apply \cite[Theorem 2.1]{psz2009} 
resp.\ \cite[Theorem 5.3.1]{moving_interfaces} to show that every 
stationary solution $(V,p_0)$ with $V \in B_{\alpha, \beta}$ is exponentially 
stable in the following sense: 
%If $(v,p)$ is a solution to 
%(\ref{eqn:min-hyd-mod}) with initial data $v_0$ and $\|v_0 - V\|_{H^2_\pi(Q_n)}$ is sufficiently small, then $v$ is exponentially convergent to some $V_\infty \in B_{\alpha, \beta}$.
%%%%%%%%%%%%%%%%%%%%%%%%%%%%%%%%%%%%%%%%%%%%%%%%%%%%%%%%%%%%%%%%%%%%%%%%%%%%%
    \begin{theorem}
    \label{per_nonlin_stab_ordered}
    Let $\Gamma_2, \beta > 0$, $\Gamma_0 \geq 0$, $\alpha < 0$ and $\lambda_0 \in \R$.
    Let $(V,p_0)$ with $V \in B_{\alpha, \beta}$ be a stationary state of 
    (\ref{eqn:min-hyd-mod}).
    Then $(V,p_0)$ is stable in $H^2_\pi(Q_n) \cap L^2_\sigma(Q_n)$ and there exists
    a $\delta > 0$ such that if $(v,p)$ is a solution to (\ref{eqn:min-hyd-mod})
    with initial data $v_0 \in H^2_\pi(Q_n) \cap L^2_\sigma(Q_n)$ and
    $\|v_0 - V \|_{H^2_\pi(Q_n)} < \delta$ then $v$ converges to some $V_\infty \in
    B_{\alpha, \beta}$ exponentially.
    \end{theorem}
%%%%%%%%%%%%%%%%%%%%%%%%%%%%%%%%%%%%%%%%%%%%%%%%%%%%%%%%%%%%%%%%%%%%%%%%%%%%%
    \begin{proof}
    In the notation of \cite[Theorem 2.1]{psz2009} or 
    \cite[Theorem 5.3.1]{moving_interfaces} we have
    $V = H^2_\pi(Q_n) \cap L^2_\sigma(Q_n)$, $X_0 = L^2_\sigma(Q_n)$,
    $X_1 = H^4_\pi(Q_n) \cap L^2_\sigma(Q_n)$ for the spaces,
    $\mathcal{E} = B_{\alpha,\beta}$ for the manifold, $u_* = V$ for
    the equilibrium, and
    \begin{align*}
    A\tilde{v} := A(v)\tilde{v} &:= \Gamma_2 \Delta^2 \tilde{v} - \Gamma_0 \tilde{v} + \alpha \tilde{v}
    \quad (\tilde{v} \in H^4_\pi(Q_n) \cap L^2_\sigma(Q_n)),\\
    F(v) &:= - \lambda_0 P(v \cdot \nabla)v - \beta P |v|^2 v.
    \end{align*}
    for $v \in H^2_\pi(Q_n) \cap L^2_\sigma(Q_n)$. By the structure of $A$ and $F$ (linear and semilinear respectively) it is obvious that
    \begin{align*}
    (A,F) \in C^1(H^2_\pi(Q_n) \cap L^2_\sigma(Q_n),
    \mathscr{L}(H^4_\pi(Q_n) \cap L^2_\sigma(Q_n), L^2_\sigma(Q_n)) 
    	\times L^2_\sigma(Q_n))
	\end{align*}      
    and to see that $A_o$ is the linearized operator of (\ref{eqn:min-hyd-mod}) at $V$.
    From Proposition~\ref{per_op_maxreg} we know that $A_o$ (hence also
    $A$) enjoys maximal $L^p$-regularity on
    $(0,T)$ for $T < \infty$.
    
    We will show that near $V$ the set of equilibria 
    $B_{\alpha,\beta}$ is a $C^1$-manifold in 
    $H^4_\pi(Q_n) \cap L^2_\sigma(Q_n)$ of dimension $n-1 \in \N$ and that
    the tangent space for $B_{\alpha,\beta}$ at $V$ equals $N(A_o)$.
    It is canonical to
    define a $C^1$-function which maps into $B_{\alpha, \beta}$:
    
    If $n = 3$ then $V \in B_{\alpha, \beta}$ can be written as
        \begin{align*}
        V = \sqrt{- \frac{\alpha}{\beta}}
            \begin{pmatrix}
            \sin(\theta) \cos(\phi)\\
            \sin(\theta) \sin(\phi)\\
            \cos(\theta)
            \end{pmatrix}
        \end{align*}
    for fixed $\theta \in [0,\pi]$ and $\phi \in [0, 2\pi)$. We define
	    the corresponding $C^1$ map
    as
        \begin{align*}
        &\Psi : [0,\pi] \times [0,2\pi) \to H^4_\pi(Q_n) \cap L^2_\sigma(Q_n),\\
        \begin{pmatrix}
        y\\
        z
        \end{pmatrix}
        &\mapsto 
        \Psi(y,z) :=
        \sqrt{- \frac{\alpha}{\beta}}
            \begin{pmatrix}
            \sin(\theta+y) \cos(\phi+z)\\
            \sin(\theta+y) \sin(\phi+z)\\
            \cos(\theta+y)
            \end{pmatrix}.
        \end{align*}
    Hence $\Psi(y,z)\in B_{\alpha,\beta}$ is a constant function in 
    $H^4_\pi(Q_n) \cap L^2_\sigma(Q_n)$ 
    for every $(y,z) \in [0,\pi] \times [0, 2\pi)$ satisfying 
	    $\Psi(0,0) = V$.
    The corresponding tangent space of $B_{\alpha,\beta}$ at $V$ is two
    dimensional and obviously given as
        \begin{align*}
        T_V B_{\alpha,\beta}  = \langle V \rangle^\perp.
        \end{align*}
    This results in 
    \begin{align*}
    	N(A_o) &= \{u \in H^4_\pi(Q_n) \cap L^2_\sigma(Q_n) : 
	u \text{ is constant and } u \perp V \}\\ 
	&= \langle V \rangle^\perp
	= T_V B_{\alpha,\beta}.
\end{align*}
    The case $n = 2$ can be proved analogously.
    
    Combining this with Lemma~\ref{per_nonlinstab_prop_a0} we proved
    that $V$ is normally stable. By \cite[Theorem 2.1]{psz2009} or
    \cite[Theorem 5.3.1]{moving_interfaces} the assertion follows.
    \end{proof}
%%%%%%%%%%%%%%%%%%%%%%%%%%%%%%%%%%%%%%%%%%%%%%%%%%%%%%%%%%%%%%%%%%%%%%%%%%%%%

%%%%%%%%%%%%%%%%%%%%%%%%%%%%%%%%%%%%%%%%%%%%%%%%%%%%%%%%%%%%%%%%%%%%%%%%%%%%%
\subsection{Ordered polar states: normal
hyperbolicity}\label{subsec_nt_psnh}
%%%%%%%%%%%%%%%%%%%%%%%%%%%%%%%%%%%%%%%%%%%%%%%%%%%%%%%%%%%%%%%%%%%%%%%%%%%%%
%%%%%%%%%%%%%%%%%%%%%%%%%%%%%%%%%%%%%%%%%%%%%%%%%%%%%%%%%%%%%%%%%%%%%%%%%%%%%

In this last subsection,
we will show for the unstable regime that the ordered polar states
are normally hyperbolic. This gives instability in the following sense: For each sufficiently small $\rho > 0$ there exists $0 < \delta \leq \rho$ such that the unique solution $v$ of (\ref{eqn:min-hyd-mod}) with initial value $v_0 \in B_{H^2}(V,\delta)$ either satisfies
	\begin{enumerate}
	\item[(i)] $dist_{H^2}(v(t_0), B_{\alpha,\beta}) > \rho$ for a finite time $t_0 > 0$ or
	\item[(ii)] $v(t)$ exists on $\R_+$ and converges at exponential rate to some $v_\infty \in B_{\alpha,\beta}$ in $H^2_\pi(Q_n) \cap L^2_\sigma(Q_n)$ as $t \to \infty$.
	\end{enumerate}
In order to prove this, we will apply the principle of
normally hyperbolic equilibria \cite[Theorem 6.1]{psz2009} resp.\ 
\cite[Theorem 5.5.1]{moving_interfaces}.
%%%%%%%%%%%%%%%%%%%%%%%%%%%%%%%%%%%%%%%%%%%%%%%%%%%%%%%%%%%%%%%%%%%%%%%%%%%%%
    \begin{theorem}
    \label{per_nonlin_hyperbolic_ordered}
    Let $\Gamma_2, \beta > 0$ and $\alpha < 0$ and $\lambda_0 \in \R$.
	The ordered polar state (\ref{orderedstate}) is normally hyperbolic if
    \begin{align}\label{cond_norm_hyper}
    	\Gamma_2 \left(\frac{2\pi}{L}\right)^4 |\ell|^4 + \Gamma_0 \left(\frac{2\pi}{L}\right)^2 |\ell|^2 \notin [2\alpha, 0], \;\; \ell \in \Z^n\setminus\{0\}
    \end{align}
    for $\Gamma_0 < 0$ and if there exists some $\ell_0 \in \Z^n$
    such that (\ref{per_linstab_ordered_eq}) holds.
    Thus, the ordered polar state (\ref{orderedstate}) is unstable 
    in $H^2_\pi(Q_n) \cap L^2_\sigma(Q_n)$ in the sense given above.    
    \end{theorem}
%%%%%%%%%%%%%%%%%%%%%%%%%%%%%%%%%%%%%%%%%%%%%%%%%%%%%%%%%%%%%%%%%%%%%%%%%%%%%
    \begin{proof}
    In order to apply \cite[Theorem 6.1]{psz2009} or 
    \cite[Theorem 5.5.1]{moving_interfaces} we have to show that $V$ is normally hyperbolic. In the proof of Theorem~\ref{per_nonlin_stab_ordered} we already showed that the ordered polar state forms a $C^1$-manifold of equilibria. Next, we characterize $N(A_o)$. Let $u \in N(A_o)$. Then
    \begin{align*}
    	(A_o u, A_o u)_{2,\pi} = \sum_{\ell \in \Z^n} |\sigma_{A_o}(\ell) \hat{u}(\ell)|^2 = 0.
    \end{align*}
    This yields $\sigma_{A_o}(\ell)\hat{u}(\ell) = 0$ for every $\ell \in \Z^n$, hence
    \begin{align*}
    	0 &=\Re\left(\overline{\hat{u}(\ell)}^T \sigma_{A_o}(\ell)\hat{u}(\ell)\right)\\
    	&=\Gamma_2 \left(\frac{2\pi}{L}\right)^4 |\ell|^4 |\hat{u}(\ell)|^2
    	+ \Gamma_0 \left(\frac{2\pi}{L}\right)^2 |\ell|^2 |\hat{u}(\ell)|^2\\
    	&+ 2 \beta \overline{\hat{u}(\ell)}^T \sigma_P(\ell) VV^T \sigma_P(\ell) \hat{u}(\ell)
    \end{align*}
    for all $\ell \in \Z^n$.
    We exploit that $\sigma_P(\ell)$ is symmetric  and $\sigma_P(\ell)
    \hat{u}(\ell) = \hat{u}(\ell)$ to obtain
    \begin{align*}
    	\Gamma_2 \left(\frac{2\pi}{L}\right)^4 |\ell|^4 |\hat{u}(\ell)|^2
    	+ \Gamma_0 \left(\frac{2\pi}{L}\right)^2 |\ell|^2 |\hat{u}(\ell)|^2
    	+ 2 \beta |V \cdot\hat{u}(\ell)|^2= 0.
    \end{align*}
    Setting $\ell = 0$ yields $V \perp \hat{u}(0)$ immediately. Moreover, for $\ell \neq 0$ and $\hat{u}(\ell) \neq 0$ we obtain
    \begin{align}\label{nonlinear-hyperbolicity-condition}
    	\Gamma_2 \left(\frac{2\pi}{L}\right)^4 |\ell|^4 
    	+ \Gamma_0 \left(\frac{2\pi}{L}\right)^2 |\ell|^2 = - \frac{2 \beta |V\cdot \hat{u}(\ell)|^2}{|\hat{u}(\ell)|^2} \in [2\alpha, 0]
    \end{align}
    since $|V|^2 = - \alpha/\beta$. Due to (\ref{cond_norm_hyper}) we have $\hat{u}(\ell) = 0$. This yields
    \begin{align*}
    	N(A_o) = \{ u \in H^4_\pi(Q_n) \cap L^2_\sigma(Q_n) : u \:\text{constant and} \: u \perp V \}
    \end{align*}
    with dimension $n-1$, such that $T_V B_{\alpha,\beta} = N(A_o)$. The
    fact that $\lambda = 0$ is semisimple follows analogously to the
    proof of Theorem \ref{per_nonlin_stab_ordered}. Finally, we have to
    verify condition (iv)'. To this end, let $\lambda = ir \in
    \sigma({A_o})$ for $r \neq 0$ and $u\neq 0$ be a corresponding eigenfunction. Testing $(ir - A_o)u$ with itself we obtain
    \begin{align*}
    \overline{\hat{u}(0)}^T(ir - 2\beta VV^T)\hat{u}(0) = 0
    \end{align*}
such that
	\begin{align*}
    \Im \left ( \overline{\hat{u}(0)}^T (ir - 2\beta VV^T)\hat{u}(0) \right )
    &= \Im \left ( ir |\hat{u}(0)|^2 - 2 \beta |V \cdot \hat{u}(0)|^2 \right) \\
    &= r|\hat{u}(0)|^2 = 0
    \end{align*}
    for $\ell = 0$ and
    \begin{align*}
    \Re \left (\overline{\hat{u}(\ell)}^T(ir - \sigma_{A_o}(\ell))\hat{u}(\ell) \right)
    &= -\Re \left (\overline{\hat{u}(\ell)}^T\sigma_{A_o}(\ell)\hat{u}(\ell) \right)
    =0
    \end{align*}
    for $\ell \neq 0$. This implies $r = 0$ if $\hat{u}(0) \neq 0$ and
	    $\hat{u}(\ell) = 0$ for $\ell \neq 0$ again by assumption
	    (\ref{cond_norm_hyper}). Consequently,
	    $\lambda = 0$. By
	    \cite[Theorem~6.1]{psz2009} or 
	    \cite[Theorem~5.5.1]{moving_interfaces} 
	    the result follows.
    \end{proof}
%%%%%%%%%%%%%%%%%%%%%%%%%%%%%%%%%%%%%%%%%%%%%%%%%%%%%%%%%%%%%%%%%%%%%%%%%%%%%
	\begin{remark}
	It is easily checked that, e.g., by setting $L = 2\pi$, $\Gamma_2 = 4$, 
	$\Gamma_0 = -5$ and $\alpha = -1/4$ all conditions of 
	Theorem~\ref{per_nonlin_hyperbolic_ordered} are satisfied, which
	yields unstable equilibria on $B_{\alpha,\beta}$. Hence, 
	the condition (\ref{cond_norm_hyper}) is meaningful.
	\end{remark}
%%%%%%%%%%%%%%%%%%%%%%%%%%%%%%%%%%%%%%%%%%%%%%%%%%%%%%%%%%%%%%%%%%%%%%%%%%%%%
\begin{remark}
        Note that a normal hyperbolic equilibrium implies the existence
	of a stable and of an unstable foliation near $V$. In fact, if $V$ is normally hyperbolic, then there exists $r > 0$ and a manifold $\mathcal{M}^s$, called the stable foliation, such that for each $v_0 \in B_{H^2}(V, r)$ we have that $v_0 \in \mathcal{M}^s$, if and only if the solution $v(v_0, t)$ exists on $\R_+$ and converges to some $W \in B_{\alpha, \beta}$ at an exponential rate. Furthermore, the projection onto the stable part of $A_o$ is exactly the projection onto the tangent space of $\mathcal{M}^s$ at $V$ (cf. \cite[Thm 3.1]{psw2013}).
        Analogously, there exists an unstable foliation $\mathcal{M}^u$ (cf. \cite[Thm. 4.1]{psw2013}). 
        \end{remark}

%%%%%%%%%%%%%%%%%%%%%%%%%%%%%%%%%%%%%%%%%%%%%%%%%%%%%%%%%%%%%%%%%%%%%%%%%%%%%
\section{Conclusion}
%%%%%%%%%%%%%%%%%%%%%%%%%%%%%%%%%%%%%%%%%%%%%%%%%%%%%%%%%%%%%%%%%%%%%%%%%%%%%
%%%%%%%%%%%%%%%%%%%%%%%%%%%%%%%%%%%%%%%%%%%%%%%%%%%%%%%%%%%%%%%%%%%%%%%%%%%%%
In this note stability resp.\ instability for the active fluid model
\eqref{eqn:min-hyd-mod}
in the periodic setting is considered. 
Depending on the values of the involved parameters 
\begin{enumerate}
\item stability resp.\ instability for the disordered state and
\item normal stability resp.\ hyperbolicity for the manifold of 
ordered polar states
\end{enumerate}
are proved.
This in particular includes instability for the ordered polar states
caused by self-propulsion, often referred to as active turbulence
and observed in many applications, see e.g.\ 
\cite{Wensink-et-al:Meso-scale-turbulence,Thampi,Doostmohammadi,Dogic}.

The observed turbulence indicates existence of an attractor,
cf.\ \cite{Wensink-et-al:Meso-scale-turbulence}. To prove this
rigorously is left as a future challenge.  

\bigskip

{\bf Acknowledgements.} \ This work is supported by the DFG (German
Science Foundation) Grant SA 1043/3-1.

%%%%%%%%%%%%%%%%%%%%%%%%%%%%%%%%%%%%%%%%%%%%%%%%%%%%%%%%%%%%%%%%%%%%%%%
%\bibliography{living_fluids}
%\bibliographystyle{plain}
%%%%%%%%%%%%%%%%%%%%%%%%%%%%%%%%%%%%%%%%%%%%%%%%%%%%%%%%%%%%%%%%%%%%%%%

\end{document}